\documentclass{amsart}
\usepackage{amsfonts,amssymb,amsmath,amsthm}
\usepackage{url}
\usepackage{enumerate}

\newtheorem{thrm}{Theorem}[section]
\newtheorem{lem}[thrm]{Lemma}

\theoremstyle{definition}

\numberwithin{equation}{section}

\author{Adrian \L ydka}
\address{Faculty of Mathematics and Computer Science\\
Adam Mickiewicz University \\
ul. Umultowska 87, 61-614 Pozna{\'n}, POLAND}
\email{adrianl@amu.edu.pl}

\begin{document}

\title[Explicit formulae for the M\"obius function of an elliptic curve]{On complex explicit formulae connected with the M\"obius function of an elliptic curve}

\begin{abstract}
We study analytic properties function $m(z, E)$, which is defined on the upper half-plane as an integral from the shifted $L$-function of an elliptic curve. We show that $m(z, E)$ analytically continues to a meromorphic function on the whole complex plane and satisfies certain functional equation. Moreover, we give explicit formula for $m(z, E)$ in the strip $|\Im{z}|<2\pi$.
\end{abstract}
\maketitle

\section{Introduction}
For a complex number $z$ from the upper half-plane let
\begin{equation}
m(z)=\frac{1}{2\pi i}\int_{C}\frac{e^{sz}}{\zeta(s)}ds,
\end{equation}
where $\zeta(s)$ denotes the classical Riemann zeta function, and the path of integration consists of the half-line $s=-\frac{1}{2}+it, \infty>t\geq0,$ the line segment $\left[-\frac{1}{2}, \frac{3}{2}\right]$ and the half-line $s=\frac{3}{2}+it, 0\leq t<\infty.$ This function was considered in \cite{bartz1} and \cite{kaczorowski1} where the following theorems were proved.

\begin{thrm}[Bartz \cite{bartz1}]\label{bartz}
The function $m(z)$ can be analytically continued to a meromorphic function on the whole complex plane and satisfies the following functional equation
\begin{equation}
m(z)+\overline{m(\overline{z})}=-2\sum_{n=1}^\infty\frac{\mu(n)}{n}\cos\left(\frac{2\pi}{n}e^{-z}\right).
\end{equation}
The only singularities of $m(z)$ are simple poles at the points $z=\log{n}$, where $n$ is a square-free natural number. The corresponding residues are
$$\text{Res}_{z=\log n} m(z)=-\frac{\mu(n)}{2\pi i}.$$
\end{thrm}

J. Kaczorowski in \cite{kaczorowski1} simplified the proof of this result and gave an explicit formula for $m(z)$ in the strip $|\Im{z}|<\pi.$

\begin{thrm}[Kaczorowski \cite{kaczorowski1}]\label{kaczorowski}
For $|\Im{z}|<\pi, z\neq\log{n}, \mu(n)\neq0$ we have
\begin{equation}
\begin{split}
m(z)&=-\sum_{n=1}^\infty\frac{\mu(n)}{n}e\left(-\frac{1}{ne^z}\right)-\frac{e^z}{2\pi i}m_0(z)\\
&\quad-\frac{1}{2i}(m_1(z)+\overline{m_1}(z))+\frac{1}{2i}(F_m(z)+\overline{F_m}(z)),
\end{split}
\end{equation}
where
\begin{align*}
m_0(z)&=\sum_{n=1}^\infty\frac{\mu(n)}{n}\frac{1}{z-\log{n}}\\
\intertext{is meromorphic on $\mathbb{C}$ and}\\
m_1(z)&=\frac{1}{2\pi i}\int\limits_C\left(\tan{\frac{\pi s}{2}}-i\right)\frac{e^{sz}}{\zeta(z)}ds,\\
F_m(z)&=\frac{1}{2\pi i}\int\limits_{1}^{1+i\infty}\left(\tan{\frac{\pi s}{2}}-i\right)\frac{e^{sz}}{\zeta(z)}ds
\end{align*}
are holomorphic in the half-plane $\Im{z}>-\pi.$
\end{thrm}
In this paper we prove analogous results for the M\"obius function of an elliptic curve over $\mathbb{Q}$ defined by the Weierstrass equation

$$E/\mathbb{Q}: y^2=x^3+ax+b, \quad a, b\in\mathbb{Q}.$$
Let $L(s, E)$ denote the $L$-function of $E$ (see for instance \cite{iwaniec1}, p.365-366). For $\sigma=\Re{s}>3/2$ we have

\begin{equation}\label{Lfunkcja}
L(s,E)=\prod_{p|N}\left(1-a_pp^{-s}\right)^{-1}\prod_{p\nmid N}\left(1-a_pp^{-s}+p^{1-2s}\right)^{-1},
\end{equation}
where $N$ is the conductor of $E.$ It is well-known that coefficients $a_p$ are real and for $p\not| N$ one has

$$a_p=p+1-\#E(\mathbb{F}_p),$$
\noindent
where $\#E(\mathbb{F}_p)$ denotes the number of points on $E$ modulo $p$ including the point at infinity, and $a_p\in\{-1, 0, 1 \},$ when $p|N$ (for details see \cite{iwaniec1}, p.365).
The M\"obius function of $E$ is defined as the sequence of the Dirichlet coefficients of the inverse of the shifted $L(s, E)$:

$$\frac{1}{L\left(s+\frac{1}{2}, E\right)}=\sum_{n=1}^\infty\frac{\mu_E(n)}{n^s}, \quad \sigma>1.$$

Using (\ref{Lfunkcja}) and the well-known Hasse inequality (see \cite{iwaniec1}, p.366, (14.32)) we easily show that $\mu_E$ is a multiplicative function satisfying Ramanujan's condition ($\mu_E(n)\ll n^\epsilon$ for every $\epsilon>0$), and moreover
$$\mu_E(p^k)=
\begin{cases}
-\frac{a_p}{\sqrt{p}}, &k=1 \\
1, &k=2 \text{ and } p\not|\Delta\\
0, &k\geq3 \text{ or } k=2 \text{ and } p|\Delta
\end{cases}$$
for every prime $p$ and positive integer $k$.

Furthermore, C. Breuil , B.Conrad, F.Diamond and R.Taylor, using the method pioneered by A. Wiles, proved in \cite{breuil1} that every $L-$function of an elliptic curve analytically continues to an entire function and satisfies the following functional equation

\begin{equation}\label{rownaniefunk}
\left(\frac{\sqrt{N}}{2\pi}\right)^s\Gamma(s)L(s, E)=\eta\left(\frac{\sqrt{N}}{2\pi}\right)^{2-s}\Gamma(2-s)L(2-s, E),
\end{equation}
where $\eta=\pm1$ is called the root number.

\noindent
In analogy to $m(z)$ we define $m(z, E)$ as follows:
$$m(z, E)=\frac{1}{2\pi i}\int_C\frac{1}{L\left(s+\frac{1}{2}, E\right)}e^{sz}ds,$$
where the path of integration consists of the half-line $s=-\frac{1}{4}+it, \infty>t\geq0$, the simple and smooth curve $l$ (which is parametrized by $\tau:[0,1]\rightarrow\mathbb{C}$ such that $\tau(0)=-\frac{1}{4}, \tau(1)=\frac{3}{2}, \Im{\tau(t)}>0$ for $t\in(0,1)$ and $F(s)$ has no zeros on $l$ and between $l$ and the real axis) and the half-line $s=\frac{3}{2}+it, 0\leq t<\infty.$

\noindent
Using (\ref{rownaniefunk}) and the Stirling's formula (see \cite{iwaniec1}, p.151, (5.112)) it is easy to see that $m(z, E)$ is holomorphic on the upper half-plane.

Our main goal in this paper is to prove the following results which are extensions of Theorems \ref{bartz} and \ref{kaczorowski}.

\begin{thrm}\label{maintheorem1}
The function $m(z, E)$ can be continued analytically to a meromorphic function on the whole complex plane and satisfies the following functional equation
\begin{equation}\label{mrownaniefunk}    m(z, E)+\overline{m}(z, E)=-\frac{2\pi}{\eta\sqrt{N}}\sum_{n=1}^{\infty}\frac{\mu_E(n)}{n}J_1\left(\frac{4\pi}{\sqrt{Nn}}e^{-\frac{z}{2}}\right)-R(z),
\end{equation}

\noindent
where $R(z)=\sum \text{Res}_{s=\beta}\frac{e^{sz}}{L\left(s+\frac{1}{2}, E\right)}$ (summation is over real zeros of $L\left(s+\frac{1}{2}, E\right)$ in $(0,1)$, if there are any) and $J_1(z)$ denotes the Bessel function of the first kind
        $$J_1(z)=\sum_{k=1}^\infty\frac{(-1)^k(z/2)^{2k+1}}{k!\Gamma(k+2)}.$$

\noindent
The only singularities of $m(z, E)$ are simple poles at the points $z=\log{n}, \mu_E(n)\neq0$ with the corresponding residues
$$\text{Res}_{z=\log{n}}m(z, E)=-\frac{\mu_E(n)}{2\pi i}.$$
\end{thrm}

Let $Y_1(z)$ be the Bessel function of the second kind and let
\begin{equation}\label{hankel}
H_1^{(2)}(z)=J_1(z)-iY_1(z)
\end{equation}
denote the classical Hankel function  (see \cite{bateman2}, p.4). Moreover, let
\begin{align*}
R^{*}(z)&=\text{Res}_{s=\frac{1}{2}}\left(\tan{\pi s}\frac{e^{sz}}{L\left(s+\frac{1}{2}, E\right)}\right)+\sum\text{Res}\left(\tan{\pi s}\frac{e^{sz}}{L\left(s+\frac{1}{2}, E\right)}\right),\\
\intertext{summation is over real zeros of $L\left(s+\frac{1}{2}, E\right)$ in $(0,1)\setminus\left\{\frac{1}{2}\right\}$ (if there are any),}\\
m_0(z, E)&=\sum_{n=1}^\infty\frac{\mu_E(n)}{n^{\frac{3}{2}}}\frac{1}{z-\log n},\\
m_1(z, E)&=\frac{1}{2\pi i}\int_\textbf{\itshape{C}}\left(\tan{\pi s}-i\right)\frac{e^{sz}}{L\left(s+\frac{1}{2}, E\right)}ds,\\
H(z, E)&=\frac{1}{2\pi i}\int_\frac{3}{2}^{\frac{3}{2}+i\infty}\left(\tan{\pi s}-i\right)\frac{e^{sz}}{L\left(s+\frac{1}{2}, E\right)}ds.\\
\end{align*}
It is plain to see that $R(z)$ and $R^{*}(z)$ are entire functions, $m_0(z, E)$ is meromorphic on the whole plane, whereas $m_1(z, E)$ and $H(z, E)$ are holomorphic for $\Im{z}>-2\pi$. With this notation we have the following result.
\begin{thrm}\label{maintheorem2}
For $z=x+iy, |y|<2\pi, x\in\mathbb{R}, z\neq\log n, \mu_E(n)\neq0$ we have
\begin{equation}\label{formula}
\begin{split}
m(z, E)&=\frac{-\pi}{\eta\sqrt{N}}\sum_{n=1}^\infty\frac{\mu_E(n)}{n}\left(H_1^{(2)}\left(\frac{4\pi}{\sqrt{Nn}}e^{-\frac{z}{2}}\right)-\frac{2}{\pi}i\left(\frac{4\pi}{\sqrt{Nn}}e^{-\frac{z}{2}}\right)^{-1}\right)\\
&-\frac{1}{2}(R(z)-iR^{*}(z))+\frac{1}{2i}\left(H(z, E)+\overline{H}(z, E)\right)-\frac{e^{\frac{3}{2}z}}{2\pi i}m_0(z, E)\\
&-\frac{1}{2i}\left(m_1(z, E)+\overline{m_1}(z, E)\right).
\end{split}
\end{equation}
\end{thrm}

\section{An auxiliary lemma}

 We need the following technical lemma.
   \begin{lem}\label{oszacowanie}
Let $z=x+iy, y>0, s=Re^{i\theta}, R\sin{\theta}\geq1,\frac{\pi}{2}\leq\theta\leq\pi.$ Then for $R\geq R(x, y)$ we have
$$\left|\frac{e^{sz}}{L\left(s+\frac{1}{2}, E\right)}\right|\leq e^{-y\frac{R}{2}}.$$
\begin{proof}
Using (\ref{rownaniefunk}), the Stirling's formula and estimate \newline $\log{L\left(\sigma+it, E\right)}\ll\log(|t|+2),\quad |\sigma|\geq\frac{3}{2},\quad|t|\geq1$ (see \cite{perelli2}, p. 304) we obtain
\begin{equation}\label{rownosc}
\log\left|\frac{e^{sz}}{L\left(s+\frac{1}{2}, E\right)}\right|=\Re \log\frac{e^{sz}}{L\left(s+\frac{1}{2}, E\right)}=2R\log{R}\cos{\theta}+Rf(\theta, x, y)+O(\log{R}),
\end{equation}
where $f(\theta, x, y)=\left(x+2\log\frac{\sqrt{N}}{2\pi}-2\right)\cos\theta-(y+2\theta-\pi)\sin\theta.$

\noindent
For $\frac{\pi}{2}\leq\theta\leq\frac{\pi}{2}+\frac{1}{\sqrt{\log R}}$ we have

$$f(\theta, x, y)=-(y+2\theta-\pi)+O\left(\frac{1}{\sqrt{\log R}}\right)$$

\noindent
and hence

$$\log\left|\frac{e^{sz}}{L\left(s+\frac{1}{2}, E\right)}\right|\leq-\frac{yR}{2}.$$
For $\frac{\pi}{2}+\frac{1}{\sqrt{\log R}}\leq\theta\leq\pi$ we have

$$|\cos{\theta}|\gg\frac{1}{\sqrt{\log R}}$$
and consequently

$$\log\left|\frac{e^{sz}}{L\left(s+\frac{1}{2}, E\right)}\right|=-2|\cos{\theta}|R\log{R}+O(R)\leq-yR\leq-\frac{yR}{2}$$
for sufficently large $R$, and the lemma easily follows.

\end{proof}
\end{lem}

\section{Proof of Theorem \ref{maintheorem1}}

       We shall first prove that $m(z, E)$ has meromorphic continuation to the whole complex plane.

\noindent
Let us write
\begin{equation}
\begin{split}\label{trzycalki}
2\pi i m(z, E)&=\int\limits_{-\frac{1}{4}+i\infty}^{-\frac{1}{4}}\frac{e^{sz}}{L\left(s+\frac{1}{2}, E\right)}ds+\int\limits_{l}\frac{e^{sz}}{L\left(s+\frac{1}{2}, E\right)}ds+\int\limits_{\frac{3}{2}}^{\frac{3}{2}+i\infty}\frac{e^{sz}}{L\left(s+\frac{1}{2}, E\right)}ds\\
&=n_1(z)+n_2(z)+n_3(z),\\
\end{split}
\end{equation}
say.

\noindent
Notice that $n_2(z)$ is an entire function.

\noindent
We compute $n_3(z)$ explicitly. Term by term integration gives
\begin{equation}
n_3(z)=-e^{\frac{3}{2}z}\sum_{n=1}^\infty\frac{\mu_E(n)}{n^\frac{3}{2}(z-\log{n})}.
\end{equation}
This shows that $n_3(z)$ is meromorphic on the whole complex plane and has simple poles at the points $z=\log{n}, \mu_E(n)\neq0$ with residues

\begin{equation}\label{residua}
\text{Res}_{z=\log{n}}n_3(z)=-\mu_E(n).
\end{equation}
\noindent
Let us now consider $n_1(z)$. Let $C_1$ consist of the half-line $s=\sigma+i, -\infty<\sigma\leq-\frac{1}{4}$ and the line segment $[-\frac{1}{4}+i, -\frac{1}{4}]$. Using lemma \ref{oszacowanie} we can write

$$n_1(z)=\int\limits_{C_1}\frac{e^{sz}}{L\left(s+\frac{1}{2}, E\right)}ds.$$
\noindent
Putting $s=\sigma+i,\quad\sigma\leq0$ in (\ref{rownosc}) we obtain
$$\left|\frac{e^{(\sigma+i)z}}{L\left(\frac{1}{2}+\sigma+i, E\right)}\right|\ll e^{-c_0|\sigma|\log(|\sigma|+2)},$$
hence $n_1(z)$ is an entire function.

\noindent
Then for $z\in\mathbb{C},$ $z\neq\log{n},$ $\mu_E(n)\neq0$ we have
\begin{equation}\label{malefunk}
m(z, E)+\overline{m}(z, E)=-\frac{1}{2\pi i}\int\limits_{\overline{C_1}\cup(-C_1)}\frac{e^{sz}}{L\left(s+\frac{1}{2}, E\right)}ds-R(z)
\end{equation}
and minus before contour denotes opposite direction.

\noindent
Using the equality (\ref{rownaniefunk}) we get
\begin{equation}
\frac{1}{2\pi i}\int\limits_{\overline{C_1}\cup(-C_1)}\frac{e^{sz}}{L\left(s+\frac{1}{2}, E\right)}ds=\frac{\pi}{\eta\sqrt{N}}\sum_{n=1}^\infty\frac{\mu_E(n)}{n}\cdot\frac{1}{2\pi i}\int\limits_{\overline{C_1}\cup(-C_1)}\frac{\Gamma\left(s+\frac{1}{2}\right)}{\Gamma\left(\frac{3}{2}-s\right)}\left(\frac{Nne^z}{(2\pi)^2}\right)^s ds.
\end{equation}

\noindent
The last integrand has simple poles at $s=-\frac{1}{2}, -\frac{3}{2}, -\frac{5}{2}, \ldots$. Computing residues we obtain

$$\frac{1}{2\pi i}\int\limits_{\overline{C_1}\cup(-C_1)}\frac{\Gamma\left(s+\frac{1}{2}\right)}{\Gamma\left(\frac{3}{2}-s\right)}\left(\frac{Nne^z}{(2\pi)^2}\right)^s ds=J_1\left(\frac{4\pi}{\sqrt{Nn}}e^{-\frac{z}{2}}\right).$$

\section{Proof of Theorem \ref{maintheorem2}}

\noindent
Let us now consider the function
$$m^{*}(z, E)=\frac{1}{2\pi i}\int\limits_\textbf{\itshape{C}}\tan(\pi s)\frac{e^{sz}}{L\left(s+\frac{1}{2}, E\right)}ds.$$
Using lemma \ref{oszacowanie} we can write
\begin{equation}
m^{*}(z, E)=\frac{1}{2\pi i}\left(\int\limits_{C_1\cup l}    +\int\limits_{\frac{3}{2}}^{\frac{3}{2}+i\infty}\right)\tan{\pi s}\frac{e^{sz}}{L\left(s+\frac{1}{2}, E\right)}ds=m_a^*(z, E)+m_b^*(z, E).
\end{equation}

\noindent
Using again estimation
$$\left|\frac{e^{(\sigma+i)z}}{L\left(\frac{1}{2}+\sigma+i, E\right)}\right|\ll e^{-c_0|\sigma|\log(|\sigma|+2)},\quad\sigma\leq0$$
and $\tan(\pi(\sigma+i))\ll1$ it is easy to see that $m_a^*(z, E)$ is an entire function.

\noindent
Moreover
$$m_b^*(z, E)=H(z, E)-\frac{e^{\frac{3}{2}z}}{2\pi}m_0(z, E).$$

\noindent
This gives the meromorphic continuation of $m^*(z, E)$ to the half-plane $\Im{z}>-2\pi$ and $m^*(z, E)$ has poles at the points $\log{n}, n=1,2,3,\ldots, \mu_E(n)\neq0,$ with residues

$$\text{Res}_{s=\log{n}}m^*(z, E)=-\frac{\mu_E(n)}{2\pi}.$$

\noindent
Now we consider the function $\overline{m^*}(z, E)$ Changing $s$ to $\overline{s}$ we get

$$\overline{m^{*}}(z, E)=\frac{1}{2\pi i}\int_{-\overline{C}}\tan{\pi s}\frac{e^{sz}}{L\left(s+\frac{1}{2}, E\right)}ds, \quad \Im z<2\pi.$$

\noindent
Further we have

$$ \overline{m^{*}}(z, E)=\frac{1}{2\pi i}\int_{-\left(\overline{C_1}\cup\overline{l}\right)}\tan{\pi s}\frac{e^{sz}}{L\left(s+\frac{1}{2}, E\right)}ds+\overline{H}(z, E)-\frac{e^{\frac{3}{2}z}}{2\pi}m_0(z, E).$$

\noindent
Then for $|\Im(z)|<2\pi$ we have
\begin{equation}
m^{*}(z, E)+\overline{m^{*}}(z, E)=-J(z, E)-\frac{e^{\frac{3}{2}z}}{\pi}m_0(z, E)+H(z, E)+\overline{H}(z, E)-R^{*}(z),
\end{equation}
where
\begin{equation}
J(z, E)=\frac{1}{2\pi i}\int\limits_{\overline{C_1}\cup(-C_1)}\tan{\pi s}\frac{e^{sz}}{L\left(s+\frac{1}{2}, E\right)}ds
\end{equation}

\noindent
Using functional equation (\ref{rownaniefunk}) we get

\begin{equation*}
\begin{split}
&J(z, E)=\frac{1}{2\pi i}\int\limits_{\overline{C_1}\cup(-C_1)}\tan(\pi s)\frac{e^{sz}\Gamma\left(s+\frac{1}{2}\right)}{\eta\Gamma\left(\frac{3}{2}-s\right)L\left(\frac{3}{2}-s\right)}\left(\frac{\sqrt{N}}{2\pi}\right)^{2s-1}ds\\
&=\frac{-2\pi}{\eta\sqrt{N}}\sum_{n=1}^\infty\frac{\mu_E(n)}{n}\left(\frac{1}{2\pi i}\int\limits_{\overline{C_1}\cup(-C_1)}\frac{\Gamma\left(s+\frac{1}{2}\right)\Gamma\left(s-\frac{1}{2}\right)}{\Gamma(s)\Gamma(1-s)}\left(\frac{e^z Nn}{4\pi^2}\right)^sds\right).
\end{split}
\end{equation*}

\noindent
The last integral we can compute using inverse Mellin transform (see \cite{paris1}, p.407)

$$\frac{1}{2\pi i}\int\limits_{\overline{C_1}\cup(-C_1)}\frac{\Gamma\left(s+\frac{1}{2}\right)\Gamma\left(s-\frac{1}{2}\right)}{\Gamma(s)\Gamma(1-s)}\left(\frac{e^z Nn}{4\pi^2}\right)^s ds=-Y_1\left(\frac{4\pi}{\sqrt{Nn}}e^{-\frac{z}{2}}\right)-\frac{2}{\pi}\left(\frac{4\pi}{\sqrt{Nn}}e^{-\frac{z}{2}}\right)^{-1}.$$

\noindent
Therefore
$$J(z, E)=\frac{2\pi}{\eta\sqrt{N}}\sum_{n=1}^\infty\frac{\mu_E(n)}{n}\left(-Y_1\left(\frac{4\pi}{\sqrt{Nn}}e^{-\frac{z}{2}}\right)-\frac{2}{\pi}\left(\frac{4\pi}{\sqrt{Nn}}e^{-\frac{z}{2}}\right)^{-1}\right).$$

\noindent
For $x\in\mathbb{R}, x\neq\log{n}$ we have
\begin{equation}
\begin{split}
\Re(m^*(x, E))=\frac{\pi}{\eta\sqrt{N}}\sum_{n=1}^\infty\frac{\mu_E(n)}{n}\left(-Y_1\left(\frac{4\pi}{\sqrt{Nn}}e^{-\frac{x}{2}}\right)-\frac{2}{\pi}\left(\frac{4\pi}{\sqrt{Nn}}e^{-\frac{x}{2}}\right)^{-1}\right)\\
-\frac{e^{\frac{3}{2}x}}{2\pi}m_0(x, E)+\frac{1}{2}\left(H(x, E)+\overline{H}(x, E)\right)-\frac{1}{2}R^{*}(x).
\end{split}
\end{equation}
Obviously

$$m^*(z, E)=im(z, E)+m_1(z, E),$$
therefore we get

\begin{equation}\label{urojone}
\begin{split}
\Im(m(x, E))=&-\frac{\pi}{\eta\sqrt{N}}\sum_{n=1}^\infty\frac{\mu_E(n)}{n}\left(-Y_1\left(\frac{4\pi}{\sqrt{Nn}}e^{-\frac{x}{2}}\right)-\frac{1}{\pi}\frac{e^{\frac{x}{2}}\sqrt{Nn}}{2\pi}\right)\\
&+\frac{e^{(\frac{3}{2}x}}{2\pi}m_0(x, E)-\frac{1}{2}\left(H(x, E)+\overline{H}(x, E)\right)\\
&+\frac{1}{2}\left(m_1(x, E)+\overline{m_1}(x, E)\right)+\frac{1}{2}R^{*}(x).
\end{split}
\end{equation}
On the other hand
\begin{equation}\label{rzeczywiste}
\Re(m(x, E))=-\frac{\pi}{\eta\sqrt{N}}\sum_{n=1}^\infty\frac{\mu_E(n)}{n}J_1\left(\frac{4\pi}{\sqrt{Nn}}e^{-\frac{x}{2}}\right)-\frac{1}{2}R^(x).
\end{equation}

The equations (\ref{urojone}) and (\ref{rzeczywiste}) imply the formula for $z\in\mathbb{R},\quad z\neq\log{n},\newline \mu_E(n)\neq0,$ and by the analytic continuation, formula (\ref{formula}) is valid in the strip $|\Im{z}|<2\pi.$

\proof[Acknowledgements]
This paper is a part of my PhD thesis. I thank my thesis advisor Prof. Jerzy Kaczorowski for
suggesting the problem and helpful discussions.


\vspace{3cm}

\begin{thebibliography}{99}

\bibitem{bartz1}
  K. Bartz, \textit{On some complex explicit formulae connected with the M\"obius function, I}, Acta Arith., \textbf{57}, 1991, no.4, 283-293.

   \bibitem{bateman2}
   H. Bateman and A. Erdelyi, \textit{Higher transcendental functions}, vol. II, Mc Graw-Hill Book Company, 1953.

  \bibitem{breuil1}
  C.Breuil, B.Conrad, F.Diamond and R.Taylor, \textit{On the modularity of elliptic curves over Q}, Journal of AMS, \textbf{14} (2001), 843--939.
  
  \bibitem{iwaniec1}
  H. Iwaniec and E. Kowalski, \textit{Analytic number theory}, AMS, 2003.
  \bibitem{kaczorowski1}
 J. Kaczorowski, \textit{Results on the M\"obius function,} J. Lond. Math. Soc., \textbf{75} (2007), no.2, 509--521.
   \bibitem{paris1}
   D. Kaminski, R. B. Paris, \textit{Asymptotics and Mellin-Barnes integrals}, Cambridge University Press, 2001.
    \bibitem{perelli2}
  A. Perelli, \textit{General $L$-functions}, Ann. Mat. Pura Appl.(4), vol. \textbf{130}, 1982, 287--306.

\end{thebibliography}
\end{document}